\documentclass[11pt,twoside]{article}
\usepackage{amsthm}
\usepackage{amsmath}
\usepackage{amssymb}
\usepackage{array}
\usepackage{amsfonts,amscd}
\usepackage{exscale}
\usepackage{eucal}
\usepackage{latexsym,amsmath}
\usepackage{indentfirst}
\newcommand{\e}{\varepsilon}

\numberwithin{equation}{section}
\newtheorem{Prop}{\bf Proposition}[section]

\begin{document}
\def \b{\Box}

\begin{center}
{\large {\bf STOCHASTIC POISSON EQUATIONS ASSOCIATED TO LIE
ALGEBROIDS AND SOME REFINEMENTS OF A PRINCIPAL BUNDLE}}
\end{center}

\begin{center}
{\bf Gheorghe IVAN and Dumitru OPRI\c S}
\end{center}

\setcounter{page}{1}

\pagestyle{myheadings}

{\small {\it Dedicated to Professor doctor docent Dan I. Papuc at
his 80th anniversary}}\\

{\small {\bf Abstract}. The aim of this paper is to present the
stochastic Poisson equations associated to Lie algebroids. The
obtained results are used for determination of stochastic Poisson
equations associated to a refinement of a principal bundle having
the affine group as structurgroup and defined by the linear
group.}
{}

\section{Introduction}
\indent\indent The stochastic Poisson equations has been
introduced by J. -M. Bismut in \cite{bism} for Brownian motions.
These have extended for semimartingales in \cite{eme8}. In the
paper \cite{laza} suggest to the study of stochastic Poisson
equations on Lie algebroids, to have care in that the dual space
of the algebroid is endowed with a Poisson structure.

In this paper we give an answer of the above question and one
obtains in a canonical way the stochastic Poisson equations on Lie
algebroids. These results are used for to write the stochastic
Poisson equations associated to the principal bundles which
compose a tissue defined by the principal bundle of affine tangent
frames on a manifold and the sequence $ GA(n,{\bf R}) \supset
GL(n,{\bf R}) \supset \{e=(\delta_{j}^{i})\}$ , studied by Dan I.
Papuc in \cite{papu}(1972; MR 53 $\#$ 4058) and Dan I. Papuc and
Ion P. Popescu in \cite{papo} (1973; MR 57 $\#$ 13739). For more
details concerning the tissues and refinements of a differentiable
principal bundles defined by closed subgroups of the structure
group can be consult the paper \cite{ivop} ( Gh. Ivan and D.
Opri\c s, 2002; MR 2005 b: 55032) and the references.

The paper is structured as follows. In Section 2, some basic facts
on manifold valued semimartingale and stochastic Poisson equations
are reviewed. In Section 3 are established the stochastic Poisson
equations on a Lie algebroid. The stochastic Poisson equations
associated to a refinement of principal bundles defined by the
affine group and linear group are described in Section 4.

The study realized in this paper may be extended to other
manifolds which are equipped with Poisson structures.

Throughout this paper all the geometrical objects like, manifolds,
maps and functions always be assumed to be smooth.

\section{ Manifold valued semimartingale and stochastic Poisson equation}

\indent\indent We recall the minimal necessary backgrounds on
stochastic differential geometry (for notation, concepts and
further details see \cite{cema}, \cite{laza}).

Let $M$ be a smooth manifold of dimension $n$. A continuous $~M-$
valued stochastic process $~\Gamma $ defined on the filtered
probability space $( \Omega, {\cal F}, P, \{{\cal F}_{t}\}_{t\geq
0} ) $ is called a {\it semimartingale } if,  for any $ f\in
C^{\infty}(M) $, the process $ f\circ \Gamma $ is a real valued
semimartingale.

 Let now $V$ be a real vector space of dimension $ r$. Let $( M, \{\cdot, \cdot\}) $ be a Poisson manifold, $ X: {\bf
R}_{+}\times \Omega \to V $ a semimartingale that takes values on
 $V$ with $ X_{0}= 0 $ ( $X_{0}$ is  the initial value of $X$ ),
and $ h : M \to V^{\ast} $ is a smooth function ($ V^{\ast} $
denotes the dual of $V$ ).

Let $ \{ e^{a} | a =\overline{1,r} \}$ be a basis of $ V^{\ast}$,
and $ h\in V^{\ast}$ such that $ h = h_{a} e^{a}.$

The {\it Hamiltonian equation with stochastic component $X$, and
Hamiltonian function $h$}, is the Stratonovich differential
equation:
\begin{equation}
\delta \Gamma^{h} = H(x, \Gamma^{h})\delta X, \label{(2.1)}
\end{equation}
defined by the Statonovich operator $ H(v,z) : T_{v}V \to T_{z}M $
given by
\begin{equation}
H(v,z)u=< e^{a},u> X_{h_{a}}(z).\label{(2.2)}
\end{equation}
We will refer to $ \Gamma^{h}$ as the {\it  Hamiltonian
semimartingale associated to $h$} with initial condition
$\gamma_{0}$, (\cite{laza}).
\begin{Prop}(\cite{laza})
Let $(M,\{\cdot,\cdot\}) $ be a Poisson manifold, $~X:{\bf R}_{+}
\times \Omega \to V $ a semimartingale and $ h : M \to V^{\ast} $
a smooth function. Let $ \Gamma_{0} $ be a ${\cal F}_{0}-$
measurable random variable and $ \Gamma^{h}$ the Hamiltonian
semimartingale associated to $h$ with initial condition
$\Gamma_{0}$. Let $ \xi^{h} $ be the corresponding maximal
stopping time. Then, for any stopping time $ \tau < \xi^{h} $ the
Hamiltonian semimartingale $ \Gamma^{h}$ satisfies
\begin{equation}
f(\Gamma_{(\tau)}^{h})- f(\Gamma_{(0)}^{h}) =
\int\limits_{0}^{\tau} \{f, h_{a}\}( \Gamma^{h}) d X^{a} +
\displaystyle\frac{1}{2}\int\limits_{0}^{\tau} \{\{f, h_{a}\},
h_{b}\} d [X^{a}, X^{b}],\label{(2.3)}
\end{equation}
for all  $ f\in C^{\infty}(M).$
\end{Prop}
From (2.3) follows
\begin{equation}
x^{i}(\Gamma_{(\tau)}^{h}) - x^{i}(\Gamma_{(0)}^{h}) =
\int\limits_{0}^{\tau} \{x^{i}, h_{a}\}( \Gamma^{h}) d X^{a} +
\displaystyle\frac{1}{2}\int\limits_{0}^{\tau} \{\{x^{i}, h_{a}\},
h_{b}\} d [X^{a}, X^{b}],\label{(2.4)}
\end{equation}
for $~i=\overline{1,n}, a,b=\overline{1,r}.$

 \markboth{Gheorghe Ivan, Dumitru Opri\c s}{Stochastic Poisson
equation associated to Lie algebroids and some refinements...}

The relations (2.4) can be written in the following form:
\begin{equation}
dx^{i} =  \{x^{i}, h_{a}\} d X^{a} + \{\{x^{i}, h_{a}\}, h_{b}\} d
[X^{a}, X^{b}],~~~ i=\overline{1,n}, a,b=\overline{1,r}.
\label{(2.5)}
\end{equation}

Let $(M,\{\cdot,\cdot\})$ be a Poisson manifold and the smooth
functions $ h_{a}\in C^{\infty}(M), a= 0,1,2,...,r $. Let $ h
:M\to {\bf R}^{r+1} $ be the Hamiltonian function and consider the
semimartingale $~X:{\bf R}_{+} \times \Omega \to {\bf R}^{r+1} $
given by  $~ X(t,\omega) = (t, B_{t}^{1}(\omega),...,
B_{t}^{r}(\omega))$, where $ B^{a}, a=\overline{1,r} $ are $r-$
independent Brownian motions. L$\acute{e}$vy's characterization of
Brownian motion shows (\cite{bism}) that $[B^{a}, B^{b}]_{t}= t
\delta^{ab}$.

In this setup, the equation (2.3) reads
\begin{equation}
f(\Gamma_{(\tau)}^{h})- f(\Gamma_{(0)}^{h}) =
\int\limits_{0}^{\tau} ( \{f, h_{a}\}( \Gamma^{h}) d X^{a} +
\delta^{ab} \{\{f, h_{a}\}, h_{b}\})dt + \int\limits_{0}^{\tau}
\{f, h_{a}\} d B^{a},\label{(2.6)}
\end{equation}
\begin{equation}
dx^{i} = ( \{x^{i}, h_{0}\} + \delta^{ab} \{\{x^{i}, h_{a}\},
h_{b}\}) dt + \{ x^{i},h_{a}\} d B^{a},\label{(2.7)}
\end{equation}
for $ i=\overline{1,n},~a,b=\overline{1,r}.$

 These equations have been studied by Bismut in \cite{bism} in
the particular case in which the Poisson manifold
$(M,\{\cdot,\dot\}$ is just the symplectic Euclidean space ${\bf
R}^{2n}$ with the canonical symplectic form.
\begin{Prop}
Let $( {\bf R}^{n},\{\cdot,\cdot\}) $ be a Poisson manifold with
$\{x^{i}, x^{j}\} = \Lambda_{k}^{ij} x^{k}, $ and $ h_{a} =
\alpha_{ai}x^{i},~ a=\overline{1,r} $ with $ \alpha_{ai}\in {\bf
R} $. The equation (2.7) is given by
\begin{equation}
dx^{i} = ( \Lambda_{\ell}^{ij} \displaystyle\frac{\partial
h_{0}}{\partial x^{j}} + \delta^{ab} \alpha_{aj}
\alpha_{bk}\Lambda_{p}^{ij}\Lambda_{p\ell}^{k}) x^{\ell}dt +
\alpha_{aj}\Lambda_{\ell}^{ij} x^{\ell}  d B^{a},\label{(2.8)}
\end{equation}
for $~i, j, k, \ell, p=\overline{1,n},~a, b =\overline{1,r}.$
\end{Prop}

The equations (2.8) are called the {\it stochastic Poisson
equations} associated to Poisson manifold $( {\bf
R}^{n},\{\cdot,\cdot\})$.

Applying the relations (2.8) for the Poisson structures defined on
${\bf R}^{3}, {\bf R}^{6}, {\bf R}^{9}$ one obtains the stochastic
Poisson equations for the rigid body on $ SO(3), SO(2,1), $ heavy
top etc. (\cite{anca}).

\section{Stochastic Poisson equations associated to a Lie algebroid}

The theory of Lie algebroids has recently proved to be extremely
fruitful in tackling some problems in the context of geometric
mechanics (\cite{laza}). Recall that the dual of a Lie algebroids
admits a canonical Poisson structure and, therefore, one can
naturally consider Hamiltonian systems on them. According to the
results and the acceptance of this new formalism we shall
investigate the consequences of having stochastic processes taking
values on their duals for mechanical purposes.

A {\it Lie algebroid} $A$ over a manifold $M$ is a vector bundle $
\pi : A \to M $  together with a Lie algebra structure $
[\cdot,\cdot ] $ on the space of sections $ Sec(A) $ and a bundle
map $ b :A\rightarrow TM $ ( called {\it anchor map}) such that:

$(i)~$ the induced map  $ b: Sec(A) \to Sec(TM)={\cal X}(M) $ is a
homomorphism of Lie algebras;

$(ii)~$ for any $ a_{1}, a_{2} \in  Sec(A) $ and smooth function $
f\in C^{\infty}(M) $, the {\it Leibniz identity} holds:
\begin{equation}
[ a_{1}, f a_{2}] = f [ a_{1}, a_{2}] +  b(a_{1})(f)
a_{2}.\label{(3.1)}
\end{equation}

For a  Lie algebroid $ (E, \pi, M,  [\cdot,\cdot], b ) $, we
consider the manifold $M$  of dimension $n$ and denote the rank of
the vector bundle $A$  with $r$.  Recalling the construction of a
canonical Poisson bracket on the dual $ A^{\ast}$ of the vector
bundle $A$ (\cite{cawe}). If one fixes local coordinates $ (
x^{i}), i=\overline{1,n} $  over a trivializing neighborhood $U
\subset M$ and choose a basis of local sections $\{e_{\alpha} |
\alpha=\overline{1,r} \}$ of the vector bundle $ A $, then the
corresponding local coordinates on $ A $ are denoted by $ ( x^{i},
y^{\alpha} ), i=\overline{1,n}, \alpha=\overline{1,r} $.

The local expression of a section $a\in Sec(A) $ with to respect
the basis $ \{ e_{\alpha}\}$  is $ a = a^{\alpha}e_{\alpha},$ with
$ a^{\alpha}\in C^{\infty}(U), ~\alpha=\overline{1,r}.$ Since $
e_{\alpha}\in Sec(A),$ we have  $ b(e_{\alpha})\in {\cal X}(U) $
and $ [e_{\alpha},e_{\beta}]\in Sec(A)$. Then there exists the
functions $b_{\alpha}^{i}, C_{\alpha \beta}^{\gamma} \in
C^{\infty}(U)$ such that:
 \begin{equation}
 \left \{ \begin{array}{lcl}
b(e_{\alpha}) &=& b_{\alpha}^{i}\displaystyle\frac{\partial}{\partial x^{i}},~~~~~~~~\hbox{for}~~~ i=\overline{1,n},~ \alpha=\overline{1,r}\\[0.3cm]
[e_{\alpha}, e_{\beta}]& =& C_{\alpha \beta}^{\gamma}
e_{\gamma},~~~~~~~\hbox{for}~~~ \alpha, \beta,
\gamma=\overline{1,r}. \label{(3.2)}
\end{array}\right.
\end{equation}
The functions $ b_{\alpha}^{i}, C_{\alpha \beta}^{\gamma}\in
C^{\infty}(U) $ given by the relations $ (3.2) $  are called the
{\it structure functions } of  the Lie algebroid $ (E,[\cdot
,\cdot], b) $ with to respect the chosen local coordinates system.

The defining relations for a Lie algebroid translate into certain
partially differential equations involving its structure
functions.

One define a Poisson structure on $ A^{\ast}$ as follows. Let $\{
\xi_{\alpha}\}$ the linear coordinates on the fibers of $
A^{\ast}$ associated with the basis of local sections $
e_{\alpha}, \alpha =\overline{1,r}$. The Poisson bracket
$\{\cdot,\cdot\}$ on $ C^{\infty}(A^{\ast})$ is defined by
\begin{equation}
 \Lambda^{ij}=\{ x^{i}, x^{j}\}= 0,~~~ \Lambda_{\alpha}^{i} = \{x^{i}, \xi_{\alpha}\} = b_{\alpha}^{i},~~~ \Lambda_{\alpha, \beta}=\{\xi_{\alpha}, \xi_{\beta}\}=
 C_{\alpha \beta}^{\gamma} \xi_{\gamma},\label{(3.3)}
\end{equation}
for $ i,j=\overline{1,n},~ \alpha, \beta, \gamma=\overline{1,r}.$

One checks that this bracket is independent of the choice of local
coordinates and basis.

Let $ a\in Sec(A)$ be a section of the vector bundle $A$. Then it
defines in a natural way a function $ f_{a}: A^{\ast} \to {\bf R}
$ which is linear in the fibers and is given by
\begin{equation}
f_{a}(x, \xi) = a^{\alpha}(x) \xi_{\alpha},~~
\alpha=\overline{1,r}.\label{(3.4)}
\end{equation}
\begin{Prop}
(\cite{cawe}) The assignement $ a\mapsto f_{a} $ defines a Lie
algebra homomorphism $ ( Sec(A), [\cdot,\cdot] ) \to ( C^{\infty}(
A^{\ast}), \{\cdot,\cdot\})$. Moreover, the Hamiltonian vector
field associated with $f_{a}$ is given by
\begin{equation}
X_{f_{a}} = b_{\beta}^{i} a^{\beta}\displaystyle\frac{\partial}{
\partial x^{i}} + ( a^{\gamma} C_{\beta \gamma}^{\lambda} - b_{\beta}^{j} \displaystyle\frac{\partial a^{\lambda}}{\partial x^{j}})
\xi_{\lambda}\displaystyle\frac{\partial }{\partial
\xi_{\beta}},~~~
  i,j=\overline{1,n},~ \beta, \gamma, \lambda=\overline{1,r}.\label{(3.5)}
\end{equation}
\end{Prop}
Let be the functions $ f_{s}:A^{\ast} \to {\bf R} $ for each $
s=\overline{1,p}$, where
\begin{equation}
f_{s}(x,\xi)= a_{s}^{\alpha}(x)\xi_{\alpha},
~~~\alpha=\overline{1,r}.\label{(3.6)}
\end{equation}
Using the relations (3.3) and (3.6), from (2.7) we obtain the {\it
stochastic Poisson equations} associated to $ h: A^{\ast} \to {\bf
R}$ and $ f_{s}, ~ s=\overline{1,p}$, given by
\begin{equation}
 \left \{ \begin{array}{lcl}
dx^{i} &=& (b_{\alpha}^{i}\displaystyle\frac{\partial h}{\partial
\xi_{\alpha}} + \delta^{su} b_{\lambda}^{k} a_{s}^{\lambda}
\displaystyle\frac{\partial}{\partial x^{k}}
( b_{\beta}^{i} a_{u}^{\beta}) )dt + b_{\beta}^{i}a_{s}^{\beta} d B^{s},\\[0.4cm]
d \xi_{\alpha}& =& (b_{\alpha}^{i}\displaystyle\frac{\partial
h}{\partial x^{i}} + C_{\alpha
\beta}^{\gamma}\xi_{\gamma}\displaystyle\frac{\partial h}{\partial
\xi_{\beta}}+ \delta^{su}b_{\gamma}^{j}\displaystyle\frac{\partial
}{\partial x^{j}}( b_{\alpha}^{i} \displaystyle\frac{\partial
a_{u}^{\gamma}}{\partial x^{i}})
a_{s}^{\varepsilon}\xi_{\varepsilon} + \\[0.4cm]
&& + \delta^{su} C_{\theta \gamma}^{\varepsilon} b_{\alpha}^{i}
\displaystyle\frac{\partial a_{u}^{\theta}}{\partial x^{i}}
a_{s}^{\gamma} \xi_{\varepsilon}) dt + ( b_{\alpha}^{i}
\displaystyle\frac{\partial a_{s}^{\lambda}}{\partial x^{i}}
\xi_{\lambda} + C_{\alpha \mu}^{\gamma} a_{s}^{\mu} \xi_{\gamma})d
B^{s}. \label{(3.7)}
\end{array}\right.
\end{equation}
Let the  tangent bundle $ TM \to M$ and cotangent bundle $
T^{\ast}M\to M$. The total space of the vector bundle $ T^{\ast}M
\oplus A^{\ast}$ has the Poisson structure $\{\cdot,\cdot \}$,
defined by
\begin{equation}
\left\{ \begin{array}{lll}
 \Lambda^{ij}=\{ x^{i}, x^{j}\}= 0, &\Lambda_{j}^{i} = \{x^{i}, p_{j}\} = \delta_{j}^{i},& \Lambda_{\alpha}^{i}=\{ x^{i},
 \xi_{\alpha}\},\\[0.2cm]
 \Lambda_{ij}=\{p_{i},p_{j}\},&\Lambda_{\alpha \beta}=\{\xi_{\alpha},
 \xi_{\beta}\}=C_{\alpha \beta}^{\gamma} \xi_{\gamma},&
\Lambda_{i\alpha}=\{p_{i},\xi_{\alpha}\}. \label{(3.8)}
\end{array}\right.
\end{equation}
\begin{Prop}
The stochastic Poisson equations defined by $ h: T^{\ast}M\oplus
A^{\ast} \to {\bf R} $ and functions $ g_{s}: T^{\ast}M\oplus
A^{\ast} \to {\bf R},~s=\overline{1,p}$, given by
\begin{equation}
\left\{\begin{array}{lcl}
 g_{s}(x,p,\xi)& =& a_{s}^{\alpha}(x)
\xi_{\alpha}+ d_{s}^{i}p_{i},~ s=\overline{1,p},\\[0.2cm]
h(x,p,\xi)& =&\displaystyle\frac{1}{2} k^{ij}(x)p_{i}p_{j}+
k^{i\alpha}(x)p_{i}\xi_{\alpha} + \displaystyle\frac{1}{2}
k^{\alpha \beta}(x)\xi_{\alpha} \xi_{\beta},\label{(3.9)}
\end{array}\right.
\end{equation}
 are
\begin{equation}
 \left \{ \begin{array}{lcl}
dx^{i} &=& (k^{ij}+ b_{\alpha}^{i}k^{j\alpha}) p_{j} +
(k^{i\beta}+ b_{\alpha}^{i}k^{\alpha \beta}) p_{\beta} +
\delta^{su} ( d_{u}^{j}+ b_{\alpha}^{j} a_{u}^{\alpha})\cdot\\[0.4cm]
&&\cdot\displaystyle\frac{\partial}{\partial x^{j}}
( a_{s}^{i}+ b_{\alpha}^{i} a_{s}^{\alpha})+ ( d_{s}^{i}+ b_{\alpha}^{i}a_{s}^{\alpha}) d B^{s}(t),\\[0.4cm]
dp_{j} &=& ( - (
\displaystyle\frac{1}{2}\displaystyle\frac{\partial
k^{h\ell}}{\partial x^{j}}p_{h} p_{\ell} + \delta^{us}(
b_{\alpha}^{m} a_{u}^{\alpha}+
d_{u}^{m})(\displaystyle\frac{\partial^{2}
a_{s}^{\gamma}}{\partial x^{m}\partial x^{j}} \xi_{\gamma}+
\displaystyle\frac{\partial^{2} d_{s}^{i}}{\partial x^{m}\partial x^{j}} p_{i})-\\[0.4cm]
&&- \delta^{su}\displaystyle\frac{\partial
d_{s}^{\alpha}}{\partial x^{j}}( \displaystyle\frac{\partial
a_{u}^{\alpha}}{\partial x^{\ell}}\xi_{\alpha}+
\displaystyle\frac{\partial d_{u}^{i}}{\partial x^{\ell}} p_{i}))
dt + ( \displaystyle\frac{\partial a_{s}^{\alpha}}{\partial
x^{j}}\xi_{\alpha}+ \displaystyle\frac{\partial
d_{s}^{i}}{\partial x^{j}} p_{i})d B^{s}(t),\\[0.4cm]
d \xi_{\alpha}& =& ( -b_{\alpha}^{i}( \displaystyle\frac{1}{2}
\displaystyle\frac{\partial k^{h\ell}}{\partial
x^{i}}p_{h}p_{\ell} + \displaystyle\frac{\partial
k^{\alpha\beta}}{\partial x^{i}}\xi{\alpha} \xi_{\beta} +
\displaystyle\frac{1}{2}\displaystyle\frac{\partial k^{j\beta}
}{\partial x^{i}}
p_{j}\xi_{\beta})-\\[0.4cm]
&&-\delta^{us}(b_{\beta}^{\ell} a_{u}^{\beta} +
d_{u}^{\ell})\cdot\displaystyle\frac{\partial}{\partial x^{\ell}}(
b_{\alpha}^{i} \displaystyle\frac{\partial a_{s}^{\beta}}{\partial
x^{i}}\xi_{\beta} + b_{\alpha}^{i} \displaystyle\frac{\partial
d_{s}^{j}}{\partial x^{i}}p_{j} )+\delta^{su} b_{\alpha}^{i}(
\displaystyle\frac{\partial d_{s}^{\ell}}{\partial x^{i}}+
b_{\gamma}^{\ell}\displaystyle\frac{\partial
a_{s}^{\gamma}}{\partial x^{i}})\cdot\\[0.4cm]
&&\cdot ( \displaystyle\frac{\partial a_{u}^{\mu}}{\partial
x^{\ell}}\xi_{\mu} + \displaystyle\frac{\partial
d_{s}^{j}}{\partial x^{i}}p_{j} ))dt - b_{\alpha}^{i}(
\displaystyle\frac{\partial a_{s}^{\beta}}{\partial x^{i}}
\xi_{\beta} + \displaystyle\frac{\partial d_{s}^{j}}{\partial
x^{i}}p_{j} )d B^{s}(t). \label{(3.10)}
\end{array}\right.
\end{equation}
\end{Prop}

\section{ Stochastic Poisson equations associated to refinement of a principal bundle having the affine group as structure group}

 We start with some definitions and results of \cite{cema} that we will use
later.

Let $ \pi_{G} : P \to M $ be a left principal bundle with the Lie
group $ G $  as structure group, where $ M=P/G.$ Let ${\cal G}$
the Lie algebra of the Lie group $G.$ The associated bundle with
standard fibre ${\cal G},$ where the action of $ G $ on ${\cal G}$
is the adjoint action is called the {\it adjoint bundle} and it is
denoted by $ \widetilde{{\cal G}}^{G}=Ad_{G}(P). $ We let $~
\widetilde{\pi}_{G}: \widetilde{{\cal G}}^{G} \to M=P/G $ denote
the projection given by $
\widetilde{\pi}_{G}([q,\xi]_{G}=[q]_{G}.$

Consider now the bundle $ TM\otimes \widetilde{\cal G}^{G} \to M$
and we assume that is given  a ( principal) connection $A^{G}$ on
the principal bundle $ \pi_{G}: P~\to~ M,~$ determined by the
local functions $\{ A_{i}^{a}(x) \}$  on $M.$  Given the basis $
\{\varepsilon_{a} | a=\overline{1,p} \}$ for the Lie algebra $
{\cal G} $ having $ \{ C_{bc}^{a} \} $ as structure constants, one
obtains the local basis $ \{ \displaystyle\frac{\partial}{\partial
x^{i}}, \varepsilon_{a} \} $ for  $ Sec(TM\otimes \widetilde{\cal
G}^{G}) $  such that $ [\varepsilon_{a}, \varepsilon_{b}]=
C_{ab}^{c} \varepsilon_{c}.~$

The corresponding covariant derivative $
\widetilde{\nabla}^{A^{G}}\xi $ of a section $
\xi=\xi^{a}\varepsilon_{a} $ and  $ X\in Sec(TM) $ reads
\begin{equation}
\widetilde{\nabla}_{X}^{A^{G}}\xi~=
X^{i}(\displaystyle\frac{\partial \xi^{a}}{\partial x^{i}} +
C_{bc}^{a} A_{i}^{Gb} \xi^{c}) \varepsilon_{a}.\label{(4.1)}
\end{equation}
The curvature $~\widetilde{B}^{A^{G}}~$ of the connection $A$ is
given by
\begin{equation}
\widetilde{B}^{A^{G}}=\displaystyle\frac{1}{2}\widetilde{B}_{ij}^{Ga}
dx^{i}\wedge dx^{j}\varepsilon_{a},~~~\hbox{where}~\label{(4.2)}
\end{equation}
\begin{equation}
\widetilde{B}_{ij}^{Ga}~=~\displaystyle\frac{\partial
A_{j}^{Ga}}{\partial x^{i}} -\displaystyle\frac{\partial
A_{i}^{Ga}}{\partial x^{j}} + C_{bc}^{a}
A_{i}^{Gb}A_{j}^{Gc}.\label{(4.3)}
\end{equation}
Let $~X_{i}\oplus \overline{\xi}_{i}\in Sec(TM\oplus
\widetilde{{\cal G}}^{G}),~i=1,2~$ be given two sections. Then
\begin{equation}
[X_{1}\oplus \overline{\xi}_{1}, X_{2}\oplus \overline{\xi}_{2}] =
[X_{1},X_{2}]\oplus \widetilde{\nabla}_{X_{1}}^{A^{G}}
\xi_{2}-\widetilde{\nabla}_{X_{2}}^{A^{G}}\xi_{1}-
\widetilde{B}^{A^{G}}(X_{1},X_{2})+
[\overline{\xi}_{1},\overline{\xi}_{2}].\label{(4.4)}
\end{equation}
For $~\{~\displaystyle\frac{\partial}{\partial x^{i}} \oplus
\varepsilon_{a},~i=\overline{1,n}, ~a=\overline{1,p}~$ we have
\begin{equation}
[\displaystyle\frac{\partial}{\partial x^{i}}\oplus
\varepsilon_{a}, \displaystyle\frac{\partial}{\partial
x^{j}}\oplus \varepsilon_{b}] = (C_{cb}^{d}A_{i}^{Gc}-
C_{ca}^{d}A_{j}^{Gc}- \widetilde{B}_{ij}^{A^{G}d} +
C_{ab}^{d})\varepsilon_{d}.\label{(4.5)}
\end{equation}
Let $ (x^{i}, \dot{x}^{i}, \xi^{a}) $ the local coordinates of $
TM\oplus \widetilde{{\cal G}}^{G} $ and $ (x^{i},p_{i},\mu_{a}) $
the local coordinates of $ T^{\ast}M\oplus \widetilde{{\cal
G}}^{G^{\ast}}. $ The structure Poisson on $ T^{\ast}M\oplus
\widetilde{{\cal G}}^{G^{\ast}} $ is given by
\begin{equation}
\left\{ \begin{array}{lll}
 \{ x^{i}, x^{j}\}= 0, & \{x^{i}, p_{j}\} = \delta_{j}^{i},& \{ p_{i},
 p_{j}\}= - B_{ij}^{c}\mu_{c},\\[0.2cm]
 \{p_{i},\mu_{a}\}= - C_{ca}^{d}A_{i}^{c}\mu_{d},& \{\mu_{a},
 \mu_{b}\}=C_{a b}^{c} \mu_{c},&
\{x^{i},\mu_{a}\}=0. \label{(4.6)}
\end{array}\right.
\end{equation}
Using the method for determination of Poisson equations in the
case of Lie algebroids one obtains the following proposition.
\begin{Prop}
The stochastic Poisson equations defined by the functions $ h:
T^{\ast}M\oplus \widetilde{{\cal G}}^{G^{\ast}} \to {\bf R} $ and
$ f:T^{\ast}M\oplus \widetilde{{\cal G}}^{G^{\ast}} \to {\bf R}$
with $~f(x,p,\mu)=a^{j}(x)p_{j} + d^{a}(x)\mu_{a}~$ are
\begin{equation}
 \left \{ \begin{array}{lcl}
dx^{i} &=& (\displaystyle\frac{\partial h}{\partial x^{i}}+
\displaystyle\frac{\partial a^{i}}{\partial x^{\ell}}a^{\ell})dt +
a^{i} d B(t),\\[0.4cm]
dp_{i} &=& ( - \displaystyle\frac{\partial h}{\partial x^{i}}-
B_{ij}^{c}\mu_{c}\displaystyle\frac{\partial h}{\partial p_{j}}-
C_{ca}^{d}\mu_{a} A_{i}^{c}
\displaystyle\frac{\partial h}{\partial \mu_{a}} + \{\{p_{i},f\},f\} ) dt -\\[0.4cm]
&&- ( B_{ij}^{c}\mu_{c}a^{j} + C_{ca}^{d}\mu_{d} A_{i}^{c} d^{a} )
dB(t),\\[0.4cm]
d \mu_{a}& =& ( C_{ca}^{d} \mu_{d} A_{j}^{c}
\displaystyle\frac{\partial h}{\partial p_{j}} + C_{ab}^{c}
\mu_{c}\displaystyle\frac{\partial h}{\partial \mu_{b}}+ \{\{\mu_{a},f\},f\} ) dt +\\[0.4cm]
&& +  ( C_{ca}^{d}\mu_{d} A_{j}^{c} a^{j} + C_{ab}^{c}\mu_{c}d^{b}
) d B(t). \label{(4.7)}
\end{array}\right.
\end{equation}
\end{Prop}
Let $ \pi_{G}:P~\to~M=P/G $ the principal bundle with the
structure group $G.$  We assume that is given a sequence $ {\cal
N}_{2} = ( G \supset K\supset \{e\})~$ of closed subgroups of $G.
$ If we denote $ \eta=( P,\pi_{G},M=P/G, G),$ then the pair
$(\eta, {\cal N}_{2}) $ determines a refinement $( \eta;
\eta_{01}, \eta_{12}) $ of $ \eta $ defined by $ K,$ where
$~\eta_{01}=(P/K,\pi_{GK}, M, G/K,G/N ) $ and $ \eta_{12}=(P,
\pi_{K}, P/K, K), $ and $ N $ is the largest normal subgroup of $
G $ included in $ K $ ( see Papuc \cite{papu}, Ivan and Opri\c s
\cite{ivop}).

Let $ A^{G} $ and $ A^{K} $ two connections on $ P $ given by the
forms $ A^{G}:TP \to {\cal G},~ A^{K}:TP \to {\cal K}, $ where $
{\cal G} $ resp.,  $ {\cal K} $  is the Lie algebra of $ G$ resp.,
$K. $

Let the adjoint bundles $ \widetilde{{\cal G}}^{G}=Ad_{G}(P) $ and
$ \widetilde{{\cal K}}^{K}=Ad_{K}(P).$ The vector bundles $
TM\oplus \widetilde{{\cal G}}^{G} \to M $ and $ T(P/K)\oplus
\widetilde{{\cal K}}^{K} \to P/K $ are called the {\it reduced
bundles associated to refinement defined by the pair} $ (\eta,
{\cal N}_{2}).$

Let us we apply the above considerations in the case when the
group $ G = GA(n, {\bf R}) $ is the affine group and $ K =
GL(n,{\bf R})$ is the linear group. We obtain thus the sequence
${\cal N}_{2}=( G = GA(n, {\bf R})\supset K= GL(n,{\bf R})\supset
\{e\}).$ The Lie algebra $ {\cal G} $ of $G $ has the base $~\{
e_{j}^{i}, e_{j}~\} $ and we have
$~[e_{j}^{i},e_{k}^{\ell}]=\delta_{k}^{i} e_{j}^{\ell}-
\delta_{j}^{\ell}
e_{k}^{i},~~[e_{j}^{i},e_{k}]=\delta_{k}^{i}e_{j},~~~[e_{i},e_{j}]=0.$
The Lie algebra $ {\cal K}$ of $ K $ has the base $ \{ e_{j}^{i}\}
$ and we have $~[e_{j}^{i},e_{k}^{\ell}]=\delta_{k}^{i}
e_{j}^{\ell}- \delta_{j}^{\ell} e_{k}^{i}.$

Let $ \pi_{G} : P \to M $ the principal bundle having  the affine
group $ G $ as structure group and the local coordinates $ (x^{i},
y_{j}^{i},y^{i}) $ on $P.$ The base of sections of the vector
bundle $ \widetilde{{\cal G}}^{G} \to M $ is
$~\e_{j}^{i}=y_{j}^{h}\displaystyle\frac{\partial}{\partial
y_{i}^{h}},~ \e_{j}=y_{j}^{h}\displaystyle\frac{\partial}{\partial
y^{h}}. $

Let $A^{G}$ a connection on the principal bundle $ \pi_{G}:
 P \to M$ given by the functions
$ (A_{kr}^{h},  A_{k}^{h}) $ on $ M. $  From (4.1) follows
\begin{equation}
\left\{\begin{array}{lcl}
\widetilde{\nabla}_{\displaystyle\frac{\partial}{\partial x^{i}}}^{A^{G}} \e_{k}^{\ell} &=& ( A_{ki}^{p}\delta_{q}^{\ell} - A_{qi}^{\ell}\delta_{k}^{p})\e_{p}^{q} - A_{i}^{\ell} \e_{k}, \\[0.8cm]
\widetilde{\nabla}_{\displaystyle\frac{\partial}{\partial x^{r}}}^{A^{K}} \e_{k} &=&  A_{kr}^{i} \e_{i}, \\[0.8cm]
\widetilde{B}^{A^{G}}& =& \displaystyle\frac{1}{2}( B_{kij}^{\ell}
d x^{i}\wedge dx^{j}\otimes \e_{\ell}^{k} + B_{ij}^{\ell}
dx^{i}\wedge d x^{j}\otimes \e_{\ell}).\label{(4.8)}
\end{array}\right.
\end{equation}
Let $ (x^{i}, p_{i}, \mu_{k}^{\ell}, \mu_{\ell}) $ the local
coordinates on $ T^{\ast}M\oplus\widetilde{{\cal G}}^{G^{\ast}}.~$
The structure Poisson is given by the following relations
\begin{equation}
\left\{\begin{array}{l}
\{x^{i},x_{j}\}=0,~~~\{x^{i},\mu_{k}^{\ell}\}=
0,~~~\{x^{i},\mu_{\ell}\}=0,~~~\{\mu_{i},\mu_{j}\}=0,\\[0.2cm]
\{x^{i},p_{j}\}=\delta_{j}^{i},~~~ \{p_{i}, p_{j}\}=
-B_{kij}^{\ell}\mu_{\ell}^{k} - B_{ij}^{\ell}\mu_{\ell},\\[0.2cm]
\{p_{i},\mu_{\ell}^{k}\} = ( A_{ki}^{p} \delta_{q}^{\ell}-
A_{qi}^{\ell}\delta_{k}^{p})\mu_{p}^{q}- A_{i}^{\ell} \mu_{k},~~~
\{p_{i},\mu_{k}\}=
A_{ki}^{p}\mu_{p},\\[0.2cm]
\{\mu_{j}^{i},\mu_{k}^{\ell}\}=
\delta_{k}^{i}\mu_{j}^{\ell}- \delta_{j}^{\ell}\mu_{k}^{i},
~~~\{\mu_{k}^{i},\mu_{j}\}=\delta_{k}^{i}\mu_{j}.\label{(4.9)}
\end{array}\right.
\end{equation}
Using the method for determination of Poisson equations in the
case of Lie algebroids one obtains the following proposition.
\begin{Prop}
The stochastic Poisson equations defined by the functions\\  $ h:
T^{\ast}M\oplus \widetilde{{\cal G}}^{G^{\ast}} \to {\bf R} $ and
$ f:T^{\ast}M\oplus \widetilde{{\cal G}}^{G^{\ast}} \to {\bf R}$
with\\[-0.5cm]
 $$f(x^{i},p_{j},\mu_{k}^{\ell}, \mu_{\ell}) = a^{j}(x)p_{j}
+ d_{\ell}^{k}(x)\mu_{k}^{\ell}+
g^{\ell}(x)\mu_{\ell}~$$\\[-0.8cm]
 are
\begin{equation}
 \left \{ \begin{array}{lcl}
dx^{i} &=& (\displaystyle\frac{\partial h}{\partial p_{i}}+
\displaystyle\frac{\partial a^{i}}{\partial x^{k}}a^{k})dt +
a^{i} d B(t),\\[0.5cm]
dp_{i} &=& ( \displaystyle\frac{\partial h}{\partial x^{i}}- (
B_{kij}^{\ell}\mu_{\ell}^{k} + B_{ij}^{\ell}\mu_{\ell})
\displaystyle\frac{\partial h}{\partial p_{j}} + ((
A_{ki}^{p}\delta_{q}^{\ell}- A_{qi}^{\ell}\delta_{k}^{p})
\mu_{p}^{q} - \\[0.5cm]
&& - A_{i}^{\ell} \mu_{k})\displaystyle\frac{\partial h}{
\partial \mu_{k}^{\ell}}+ A_{ki}^{p} \mu_{p} \displaystyle\frac{\partial h}{\partial
\mu_{k}}) + \{\{p_{i},f\},f\} ) dt + \{ p_{i}, f\}
dB(t),\\[0.5cm]
d \mu_{k}^{\ell}& =& ( (( A_{ki}^{p}\delta_{q}^{\ell} -
A_{qi}^{\ell}\delta_{k}^{p})  \mu_{p}^{q} -  A_{i}^{\ell}\mu_{k} )
d_{\ell}^{k} + A_{ki}^{p} \mu_{p}^{\ell} g^{i}+\\[0.5cm]
&&+ \{\{\mu_{k}^{\ell},f\},f\} )
dt + \{\mu_{k}^{\ell},f\}d B(t),\\[0.5cm]
 d \mu_{i}& =& (
- A_{ik}^{p}\mu_{p} a^{k} - \mu_{i} \delta_{k}^{\ell} d_{\ell}^{k}
+ \{\{\mu_{i},f\},f\} ) dt + \{\mu_{i},f\} d B(t). \label{(4.10)}
\end{array}\right.
\end{equation}
\end{Prop}
Let $ \pi_{K} : P \to P/K $ the principal bundle having the affine
group $ K= GL(n,{\bf R})$  as structure group and the local
coordinates $ (x^{i}, q^{i}) $ on $ P/K. $ The base of sections of
the vector bundle $~ \widetilde{{\cal K}}^{K} \to P/K $ is
$~\varepsilon_{j}^{i}=y_{j}^{h}\displaystyle\frac{\partial}{\partial
x_{i}^{h}}.$

Let $ A^{K} $ a connection on the principal bundle $ \pi_{K}: P
\to P/K $ given by the functions $ (A_{ij}^{k}, B_{ij}^{k}) $ on
$P/K.$  From the relations (4.1) follows:
\begin{equation}
\left\{\begin{array}{cll}
\widetilde{\nabla}_{\displaystyle\frac{\partial}{\partial x^{i}}}^{A^{K}} \varepsilon_{k}^{\ell} &=& ( A_{ki}^{p}\delta_{q}^{\ell} - A_{qi}^{\ell}\delta_{k}^{p})\varepsilon_{p}^{q} \\[0.8cm]
\widetilde{\nabla}_{\displaystyle\frac{\partial}{\partial q^{i}}}^{A^{K}} \varepsilon_{k}^{\ell} &=& ( B_{ki}^{p}\delta_{q}^{\ell} - B_{qi}^{\ell}\delta_{k}^{p})\varepsilon_{p}^{q} \\[0.8cm]
\widetilde{B}^{A^{K}} &=& \displaystyle\frac{1}{2}( B_{kij}^{\ell}
d x^{i}\wedge dx^{j}+B_{kij}^{\ell} dq^{i}\wedge d q^{j} +
B_{kij}^{\ell} d x^{i}\wedge dq^{j})\otimes \varepsilon_{\ell}^{k}
.\label{(4.11)}
\end{array}\right.
\end{equation}
Let $ (x^{i}, q^{i}, \dot{x}^{i}, \dot{q}^{i}, \xi_{k}^{\ell}) $
the local coordinates on $ T(P/K)\oplus \widetilde{{\cal K}}^{K} $
and $ (x^{i}, q^{i},p_{i},\lambda_{i}, \mu_{k}^{\ell} )$ the local
coordinates on $ T^{\ast}(P/K)\oplus\widetilde{{\cal
K}}^{K^{\ast}}. $ The structure Poisson is given by the following
relations:
\begin{equation}
\left\{\begin{array}{l}
\{x^{i},x^{j}\}= 0, ~~~\{x^{i},q^{k}\}=0,~~~\{x^{i},p_{j}\}=\delta_{j}^{i},~~~\{x^{i}, \lambda_{j}\}= \delta_{j}^{i},\\[0.2cm]
\{ x^{i}, \mu_{k}^{\ell}\}=0,~~~
\{q^{i}, q^{j}\}=0,~~~\{q^{i}, p_{j}\}=0,~~~\{q^{i}, \lambda_{j}\}=0,\\[0.2cm]
\{q^{i}, \mu_{k}^{\ell}\}=0,~~~ \{p_{i}, p_{j}\}= -
\displaystyle\frac{1}{2}
B_{kij}^{\ell}\mu_{\ell}^{k},~~~\{p_{i},\lambda_{j}\}=-\displaystyle\frac{1}{2}
B_{kij}^{\ell}\mu_{\ell}^{k},\\[0.2cm]
\{p_{i},\mu_{k}^{\ell}\}=(A_{ki}^{p}\delta_{q}^{\ell}-A_{qi}^{\ell}\delta_{k}^{p})\mu_{p}^{q},~~~
\{\lambda_{i},\lambda_{j}\}= -\displaystyle\frac{1}{2}
B_{kij}^{\ell}\mu_{\ell}^{k},\\[0.2cm]
\{\lambda_{i},\mu_{k}^{\ell}\}=(B_{ki}^{p}\delta_{q}^{\ell}-
B_{qi}^{\ell}\delta_{k}^{p}) \mu_{p}^{q},~~~
\{\mu_{j}^{i},\mu_{k}^{\ell}\}= \delta_{k}^{i}\mu_{j}^{\ell}-
\delta_{j}^{\ell}\mu_{k}^{i}.\label{(4.12)}
\end{array}\right.
\end{equation}
Using the method for determination of Poisson equations in the
case of Lie algebroids one obtains the following proposition.
\begin{Prop}
The stochastic Poisson equations defined by the functions\\  $ h:
T^{\ast}(P/K)\oplus \widetilde{{\cal K}}^{K^{\ast}} \to {\bf R} $
and $ f:T^{\ast}(P/K)\oplus \widetilde{{\cal K}}^{K^{\ast}} \to
{\bf R}$
with\\[-0.5cm]
 $$f(x^{i}, q^{i}, p_{j}, \lambda_{j}, \mu_{k}^{\ell}) = a^{j}(x,q)p_{j}
+ d^{j}(x,q)\lambda_{j}+
g_{k}^{j}(x,q)\mu_{j}^{k}~$$\\[-0.8cm]
 are
\begin{equation}
 \left \{\begin{array}{lcl}
dx^{i} &=& (\displaystyle\frac{\partial h}{\partial p_{i}}+
\{\{x^{i}, f\}, f\})dt + \{x^{i}, f\}d B(t),\\[0.5cm]
dp_{i} &=& ( - \displaystyle\frac{\partial h}{\partial x^{i}}-
\displaystyle\frac{1}{2} B_{kij}^{\ell}\mu_{\ell}^{k} + (
A_{ki}^{p}\delta_{q}^{\ell}-
A_{qi}^{\ell}\delta_{k}^{p})\mu_{p}^{q}
\displaystyle\frac{\partial h}{\partial \mu_{k}^{\ell}} + \\[0.5cm]
&& + \{\{p_{i},f\},f\} ) dt + \{ p_{i}, f\} dB(t),\\[0.5cm]
dq^{i} &=& (\displaystyle\frac{\partial h}{\partial \lambda_{i}}+
\{\{q^{i}, f\}, f\})dt + \{q^{i}, f\} d B(t),\\[0.5cm]
d\lambda_{i} &=& ( - \displaystyle\frac{\partial h}{\partial
q^{i}}- \displaystyle\frac{1}{2} B_{kij}^{\ell}\mu_{\ell}^{k} + (
B_{ki}^{p}\delta_{q}^{\ell}-
B_{qi}^{\ell}\delta_{k}^{p})\mu_{p}^{q}
\displaystyle\frac{\partial h}{\partial \mu_{k}^{\ell}} + \\[0.5cm]
&& + \{\{\lambda_{i},f\},f\} ) dt + \{ \lambda_{i}, f\} dB(t),\\[0.5cm]
d \mu_{k}^{\ell}& =& ( -( A_{kj}^{p}\delta_{q}^{\ell} +
A_{qj}^{\ell}\delta_{k}^{p})  \mu_{p}^{q} \displaystyle\frac{\partial h}{\partial p_{j}} -( B_{kj}^{p} \delta_{q}^{\ell} -
B_{kj}^{\ell} \delta_{k}^{p})\mu_{p}^{q} \displaystyle\frac{\partial h}{\partial \lambda_{j}} +\\[0.5cm]
 && +(\delta_{j}^{\ell}\mu_{k}^{i} - \delta_{k}^{i} \mu_{j}^{\ell})\displaystyle\frac{\partial h}{\partial \mu_{j}^{i}}  +
 \{\{\mu_{k}^{\ell},f\},f\} ) dt + \{\mu_{k}^{\ell},f\}d
B(t). \label{(4.13)}
\end{array}\right.
\end{equation}
\end{Prop}
The study of equations (3.10), (4.7), (4.10) and (4.13) enable by
choosing of the functions $h$ and $ f_{a}$.

\vspace*{0.4cm}

\hspace*{0.7cm}West University of Timi\c soara\\
\hspace*{0.7cm} Department of Mathematics\\
\hspace*{0.7cm} Bd. V. P{\^a}rvan,no.4, 300223, Timi\c soara, Romania\\
\hspace*{0.7cm}E-mail: ivan@math.uvt.ro; miticaopris@yahoo.com

\end{document}